\documentstyle[pb-diagram]{article}
\newcommand{\bean}{\[\begin{array}{rcl}}
\newcommand{\eean}{\end{array}\]}
\newcommand{\N}{I\!\!N}
\newcommand{\R}{I\!\!R}
\newcommand{\Z}{Z\!\!\!Z}
\newcommand{\PP}{I\!\!P}
\newcommand{\C}{\,I\!\!\!\!C}

\newcommand{\veee}{\scriptscriptstyle\vee}

\newcommand{\surj}{\longrightarrow\hspace{-1.5em}\longrightarrow}

\newcommand{\ko}{\overline}
\newcommand{\ku}{\underline}

\newcommand{\kd}{\displaystyle}
\newcommand{\kf}{\footnotesize}

\newcommand{\keps}{\varepsilon}

\newcommand{\kr}{r}
\newcommand{\tM}{\tilde{M}}
\newcommand{\tN}{\tilde{N}}
\newcommand{\tsigma}{{\tilde{\sigma}}}
\newcommand{\ttau}{{\tilde{\tau}}}

\newcommand{\A}[2]{A\!\!\!\!\!A_{#2}({#1})}
\newcommand{\G}[2]{I\!\!L_{#2}({#1})}
\newcommand{\pn}{n}
\newcommand{\ksp}{G}
\newcommand{\an}{M}
\newcommand{\aj}{j}
\newcommand{\ak}{k}
\newcommand{\al}{l}

\newcommand{\di}{i}

\newcommand{\dn}{N}
\newcommand{\sm}{m}
\newcommand{\si}{i}
\newcommand{\sj}{j}
\newcommand{\sk}{k}
\newcounter{Abschnitt}[section]
\newcommand{\neu}[1]{\protect\refstepcounter{Abschnitt}\protect
   \label{#1}\vspace{1ex}
   {\bf (\protect\arabic{section}.\protect\arabic{Abschnitt})}
                     $\qquad$}
\newcommand{\zitat}[2]{(\protect\ref{#1}.\protect\ref{#1-#2})}

\begin{document}
\title{One-parameter families containing three-dimensional toric
Gorenstein singularities}

\author{Klaus Altmann
%\thanks{This paper was written at M.I.T.\ and supported by a DAAD-fellowship.}
     \\
     \small Institut f\"ur reine Mathematik, Humboldt-Universit\"at zu Berlin
     \vspace{-0.7ex}\\ \small Ziegelstr.~13A,
     D-10099 Berlin, Germany. \vspace{-0.7ex}\\ \small E-mail:
     altmann@mathematik.hu-berlin.de}
\date{}
\maketitle

%\tableofcontents
%\par
%\vspace{2ex}

%%%%%%%%%%%%
%
%  Introduction
%
%%%%%%%%%%%%%
\section{Introduction}\label{Int}

%%%%%%%%%%
% (Int.1)
%%%%%%%%

\neu{Int-1}
Let $\sigma$ be a rational, polyhedral cone. It induces a (normal) affine toric
variety $Y_\sigma$ which may have singularities. We would like to investigate
its deformation theory.
The vector space $T^1_Y$ of infinitesimal deformations is multigraded, and
its homogeneous pieces can be determined by combinatorial formulas developed
in \cite{T1}.\\
If $Y_\sigma$ only has an isolated Gorenstein singularity, then we can say
even more (cf.\ \cite{Tohoku}, \cite{T2}):
$T^1$ is concentrated in one single multidegree,
the corresponding homogeneous piece
allows an elementary geometric description in terms of Minkowski
summands of a certain lattice polytope, and it is even possible (cf.\ \cite{versal})
to obtain the entire versal deformation of $Y_\sigma$.\\
\par

%%%%%%%%%%
% (Int.2)
%%%%%%%%

\neu{Int-2}
The first aim of the present paper is to provide a geometric interpretation
of the $T^1$-formula for arbitrary toric singularities in every multidegree.
This can be done again in terms of Minkowski summands of certain polyhedra.
However, they neither need to be compact, nor do their vertices have to be
contained in the lattice anymore (cf.\ \zitat{T1}{7}).\\
In \cite{Tohoku} we have studied so-called toric deformations
only existing in negative (i.e.\ $\in -\sigma^{\veee}$) multidegrees. They are
genuine deformations with smooth parameter space, and they are characterized
by the fact that their total space is still toric.
Now, having a new description of $T^1_Y$, we will describe in Theorem \zitat{Gd}{3} the
Kodaira-Spencer map in these terms.\\
Moreover, using partial modifications of our singularity $Y_\sigma$, we
extend in \zitat{Gd}{5}
the construction of genuine deformations to non-negative degrees.
Despite the fact that the total spaces are no longer toric, we can still
describe them and their Kodaira-Spencer map combinatorially.\\
\par

%%%%%%%%%%
% (Int.3)
%%%%%%%%

\neu{Int-3}
Afterwards, we focus on three-dimensional, toric Gorenstein singularities. As
already mentioned, everything is known in the isolated case. However, as soon
as $Y_\sigma$ contains one-dimensional singularities (which then have to be of
transversal type A$_k$), the situation changes dramatically.
In general,
$T^1_Y$ is spread into infinitely many multidegrees.
Using our geometric description of the $T^1$-pieces, we detect
in \zitat{3G}{3} all
non-trivial ones and determine their dimension (which will be one in most
cases). The easiest example of that
kind is the cone over the weighted projective plane $\PP(1,2,3)$
(cf.\ \zitat{3G}{4}).\\
At least at the moment, it seems to be impossible to describe the entire versal
deformation; it is an infinite-dimensional space. However, the
infinitesimal deformations corresponding to the one-dimensional homogeneous
pieces of $T^1_Y$ are unobstructed, and we lift them in \zitat{3G}{5} to genuine
one-parameter families. Since the corresponding multidegrees are in general
non-negative, this can be done using the construction introduced in
\zitat{Gd}{5}. See section \zitat{3G}{8} for a corresponding sequel of
example $\PP(1,2,3)$.\\
Those one-parameter families form a kind of skeleton of the entire versal deformation. The
most important open questions are the following: Which of them belong
to a common irreducible component of the base space? And, how could those
families be combined to find a general fiber (a smoothing of $Y_\sigma$)
of this component? The answers to these questions would provide important
information about three-dimensional flips.\\
\par

%%%%%%%%%%%%
%
%  T1
%
%%%%%%%%%%%%%
\section{Visualizing $T^1$}\label{T1}

%%%%%%%%%%
% (T1.1)
%%%%%%%%

\neu{T1-1} {\em Notation:}
As usual when dealing with toric varieties, denote by $N$, $M$ two mutually
dual lattices (i.e.\ finitely generated, free abelian groups), by
$\langle\,,\,\rangle:N\times M\to\Z$ their perfect pairing, and by
$N_{\R}$, $M_{\R}$ the corresponding $\R$-vector spaces obtained by extension of
scalars.\\
Let $\sigma\subseteq N_{\R}$ be the polyhedral cone with apex in $0$ given
by the fundamental generators $a^1,\dots,a^\an\in N$.
They are assumed to be primitive, i.e.\ they are not proper multiples
of other elements from $N$. We will write $\sigma=\langle a^1,\dots,a^\an\rangle$.\\
The dual cone
$\sigma^{\veee}:=\{r\in M_{\R}\,|\; \langle\sigma,r\rangle\geq 0\}$
is given by the inequalities assigned to $a^1,\dots,a^\an$.
Intersecting $\sigma^{\veee}$ with the lattice $M$ yields a
finitely generated semigroup.
Denote by $E\subseteq\sigma^{\veee}\cap M$ its minimal generating set,
the so-called Hilbert basis. Then, the affine toric variety
$Y_\sigma:=\mbox{Spec}\C[\sigma^{\veee}\cap M]\subseteq\C^E$ is given by
equations assigned to the linear dependencies among elements of $E$. See
\cite{Oda} for a detailed introduction into the subject of toric varieties.\\
\par

%%%%%%%%%%
% (T1.2)
%%%%%%%%%

\neu{T1-2}
Most of the relevant rings and modules for $Y_\sigma$ are $M$-(multi)graded.
So are the modules $T^i_Y$, which are important for describing infinitesimal
deformations and their obstructions. Let $R\in M$, then in \cite{T1} and \cite{T2} we
have defined the sets
\[
E_\aj^R:=\{ r\in E\,|\; \langle a^\aj,r\rangle<\langle a^\aj,R\rangle\}
\quad
(\aj=1,\dots,\an)\,.
\]
They provide the main tool for building a complex
$\mbox{span}(E^R)_{\bullet}$ of free Abelian groups with the usual differentials
via
\[
\mbox{span}(E^R)_{-k} := \!\!\bigoplus_{\begin{array}{c}
\tau<\sigma \mbox{ face}\\ \mbox{dim}\, \tau=k \end{array}}
\!\!\!\!\!\mbox{span}(E^R_{\tau})\quad
\mbox{with} \quad
\renewcommand{\arraystretch}{1.5}
\begin{array}[t]{rcl}
E_0^R &:=& \bigcup_{\aj=1}^N E_\aj^R\; ,\mbox{ and}\\
E^R_{\tau} &:=& \bigcap_{a^\aj \in \tau} E_\aj^R \; \mbox{ for faces }
\tau < \sigma\,.
% \vspace{1ex}
\end{array}
\]

{\bf Theorem:} (cf.\ \cite{T1}, \cite{T2})
{\em
For $i=1$ and, if $Y_\sigma$ is additionally smooth in codimension two, also for $i=2$,
the homogeneous pieces of $T^i_Y$ in degree $-R$ are
\[
T^i_Y(-R)=H^i\Big(\mbox{\em span}(E^R)_\bullet^\ast\otimes_{\Z}\C\Big)\,.
\vspace{-1ex}
\]
}
\par

In particular, to obtain $T^1_Y(-R)$, we need to determine the vector spaces
$\mbox{span}_{\C}E_\aj^R$ and $\mbox{span}_{\C}E_{\aj\ak}^R$, where
$a^\aj$, $a^\ak$ span a two-dimensional face of $\sigma$. The first one
is easy to get:
\[
\mbox{span}_{\C}E_\aj^R =
\left\{
\begin{array}{ll}
0 & \mbox{if } \langle a^\aj ,R\rangle \leq 0\\
(a^\aj )^\bot & \mbox{if } \langle a^\aj ,R\rangle =1\\
M_{\C} & \mbox{if } \langle a^\aj ,R\rangle \geq 2\, .
\end{array}
\right.
\]
The latter is always contained in
$(\mbox{span}_{\C}E_\aj^R)\cap(\mbox{span}_{\C}E_\ak^R)$ with codimension
between $0$ and $2$. As we will see in the upcoming example, its actual size
reflects the infinitesimal deformations
of the two-dimensional cyclic quotient singularity assigned to
the plane cone spanned by $a^\aj$, $a^\ak$.
(These singularities are exactly the transversal types of the two-codimensional
ones of $Y_\sigma$.)\\
\par

%%%%%%%%%%
% (T1.3)
%%%%%%%%%

\neu{T1-3}
{\bf Example:}
If $Y(n,q)$ denotes the two-dimensional
quotient of $\C^2$ by the $^{\kd\Z}\!\!/_{\!\!\kd n\Z}$-action
via
$\left(\!\begin{array}{cc}\xi& 0\\ 0& \xi^q\end{array}\!\right)$
($\xi$ is a primitive $n$-th root of unity),
then $Y(n,q)$ is a toric variety and may be given by the cone $\sigma=\langle(1,0);
(-q,n)\rangle\subseteq \R^2$.
The set
$E\subseteq \sigma^{\veee}\cap \Z^2$ consists of the lattice points
$r^0,\ldots,r^w$ along the compact faces of the boundary of
$\mbox{conv}\big((\sigma^{\veee}\setminus\{0\})\cap\Z^2\big)$. There are
integers $a_v\geq 2$ such that
$r^{v-1}+r^{v+1}=a_v\,r^v$ for $v=1,\dots,w-1$. They
may be obtained by expanding
$n/(n-q)$ into a negative continued fraction
(cf.\ \cite{Oda}, \S (1.6)).\\
Assume $w\geq 2$, let $a^1=(1,0)$ and $a^2=(-q,n)$.
Then, there are only two sets $E_1^R$ and $E_2^R$ involved, and the previous theorem
states
\[
T^1_Y(-R) = \left( \left. ^{\kd (\mbox{span}_{\C}E^R_1)\cap (\mbox{span}_{\C}E^R_2)}\!
\right/ \! {\kd \mbox{span}_{\C}(E_1^R\cap E_2^R)}\right)^\ast\,.
\]
Only three different types of $R\in\Z^2$ provide a non-trivial contribution
to $T^1_Y$:
\begin{itemize}
\item[(i)] $R=r^1$ (or analogously $R=r^{w-1}):\;$
$\mbox{span}_{\C}E^R_1 =(a^1)^\bot$,
$\mbox{span}_{\C}E^R_2 = \C^2 \;(\mbox{or }(a^2)^\bot,\,\mbox{if } w=2)$,
and $\mbox{span}_{\C}E^R_{12}=0$.
Hence, $\mbox{dim}\, T^1(-R)=1$ (or $=0$, if $w=2$).
\item[(ii)] $R=r^v$ $(2\le v\le w-2)$:\quad
$\mbox{span}_{\C}E^R_1 = \mbox{span}_{\C}E^R_2 = \C^2\,$, and
$\mbox{span}_{\C}E^R_{12}=0$.
Hence, we obtain $\mbox{dim}\, T^1(-R)=2$.
\item[(iii)] $R=p\cdot r^v$  ($1\le v\le w-1$, $\;2\le p<a_v$ for $w\ge 3$;
or $v=1=w-1$, $\;2\le p\le a_1$ for $w=2$):\quad
$\mbox{span}_{\C}E^R_1 = \mbox{span}_{\C}E^R_2 = \C^2\,$, and
$\mbox{span}_{\C}E^R_{12}=\C\cdot R\,$.
In particular, $\mbox{dim}\, T^1(-R)=1$.
\vspace{1ex}
\end{itemize}

%%%%%%%%%%
% (T1.4)
%%%%%%%%%

\neu{T1-4}
{\bf Definition:}
{\em For two polyhedra $Q', Q''\subseteq \R^n$ we define their Minkowski sum
as the polyhedron
$Q'+Q'':= \{p'+p''\,|\; p'\in Q', p''\in Q''\}$. Obviously, this notion also makes
sense for translation classes of polyhedra in arbitrary affine spaces.}\\
\par

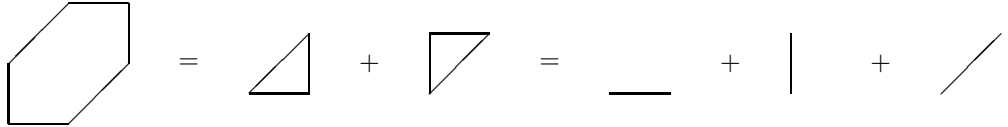
\begin{center}
\unitlength=0.40mm
\linethickness{0.4pt}
\begin{picture}(330.00,40.00)
\put(0.00,0.00){\line(1,0){20.00}}
\put(20.00,0.00){\line(1,1){20.00}}
\put(40.00,20.00){\line(0,1){20.00}}
\put(40.00,40.00){\line(-1,0){20.00}}
\put(20.00,40.00){\line(-1,-1){20.00}}
\put(0.00,20.00){\line(0,-1){20.00}}
\put(80.00,10.00){\line(1,0){20.00}}
\put(100.00,10.00){\line(0,1){20.00}}
\put(100.00,30.00){\line(-1,-1){20.00}}
\put(140.00,30.00){\line(0,-1){20.00}}
\put(140.00,10.00){\line(1,1){20.00}}
\put(160.00,30.00){\line(-1,0){20.00}}
\put(200.00,10.00){\line(1,0){20.00}}
\put(260.00,10.00){\line(0,1){20.00}}
\put(310.00,10.00){\line(1,1){20.00}}
\put(60.00,20.00){\makebox(0,0)[cc]{$=$}}
\put(120.00,20.00){\makebox(0,0)[cc]{$+$}}
\put(240.00,20.00){\makebox(0,0)[cc]{$+$}}
\put(180.00,20.00){\makebox(0,0)[cc]{$=$}}
\put(290.00,20.00){\makebox(0,0)[cc]{$+$}}
\end{picture}
\end{center}

Every polyhedron $Q$ is decomposable into the Minkowski sum
$Q=Q^{\mbox{\kf c}}+Q^{\infty}$ of a (compact) polytope $Q^{\mbox{\kf c}}$
and the so-called cone of unbounded directions $Q^{\infty}$.
The latter one is uniquely determined by $Q$,
whereas the compact summand is not. However,
we can take for $Q^{\mbox{\kf c}}$ the minimal one - given as the convex hull
of the vertices of $Q$ itself.
If $Q$ was already compact, then $Q^{\mbox{\kf c}}=Q$ and $Q^{\infty}=0$.
\vspace{1ex}\\
{\em
A polyhedron $Q'$ is called a Minkowski summand of $Q$ if there is a $Q''$ such
that $Q=Q'+Q''$ and if, additionally, $(Q')^{\infty}= Q^{\infty}$.}\\
In particular, Minkowski summands always have the same cone of unbounded directions and,
up to dilatation (the factor $0$ is allowed), the same compact edges as the
original polyhedron.\\
\par

%%%%%%%%%%
% (T1.5)
%%%%%%%%%

\neu{T1-5}
The {\em setup for the upcoming sections} is the following:
Consider the cone $\sigma\subseteq N_{\R}$ and fix some element $R\in M$.
Then $\A{R}{}:= [R=1]:= \{a\in N_{\R}\,|\; \langle a,R\rangle =1\}\subseteq N_{\R}$ is
an affine space; if $R$ is primitive, then it comes with a lattice
$\G{R}{}:= [R=1]\cap N$. The assigned vector space is $\A{R}{0}:=[R=0]$; it is always
equipped with the lattice $\G{R}{0}:= [R=0]\cap N$.
We define the cross cut of $\sigma$ in degree $R$ as the polyhedron
\[
Q(R):= \sigma\cap [R=1]\subseteq \A{R}{}\,.
\]
It has the
cone of unbounded directions $Q(R)^{\infty}=\sigma\cap \A{R}{0}\subseteq N_{\R}$.
The compact part $Q(R)^{\mbox{\kf c}}$ is given by its vertices
$\bar{a}^\aj:=a^\aj/\langle a^\aj,R\rangle$, with $\aj$
meeting $\langle a^\aj,R\rangle\geq 1$.
A trivial but nevertheless important observation is the following:
The vertex $\bar{a}^\aj $ is a lattice point (i.e.\ $\bar{a}^\aj \in \G{R}{}$),
if and only if $\langle a^\aj , R \rangle =1$.\\
Fundamental generators of $\sigma$ contained in
$R^\bot$ can still be ``seen'' as edges in $Q(R)^{\infty}$, but those with
$\langle \bullet, R\rangle <0$ are ``invisible'' in $Q(R)$. In particular, we can
recover the cone $\sigma$ from $Q(R)$ if and only if $R\in \sigma^{\veee}$.\\
\par

%%%%%%%%%%
% (T1.6)
%%%%%%%%%

\neu{T1-6}
Denote by $d^1,\dots,d^\dn\in R^\bot\subseteq N_{\R}$ the compact edges of $Q(R)$.
Similar to \cite{versal}, \S 2, we assign to each compact 2-face
$\keps<Q(R)$ its sign vector $\ku{\keps}\in \{0,\pm 1\}^\dn$ by
\[
\keps_\di := \left\{
\begin{array}{cl}
\pm 1 & \mbox{if $d^\di$ is an edge of $\keps$}\\
0     & \mbox{otherwise}
\end{array} \right.
\]
such that the oriented edges $\keps_\di\cdot d^\di$ fit into a cycle along the boundary
of $\keps$.  This determines $\ku{\keps}$ up to sign, and any choice will do.
In particular, $\sum_\di \keps_\di d^\di =0$.\\
\par

{\bf Definition:}
{\em
For each $R\in M$ we define the vector spaces
\vspace{-2ex}
\[
\renewcommand{\arraystretch}{1.5}
\begin{array}{rcl}
V(R) &:=& \{ (t_1,\dots,t_\dn)\, |\; \sum_\di t_\di  \,\keps_\di  \,d^\di  =0\;
\mbox{ for every compact 2-face } \keps <Q(R)\}\\
W(R) &:=& \R^{\#\{\mbox{\kf $Q(R)$-vertices not belonging to $N$}\}}\,.
\end{array}
\vspace{-2ex}
\]}

Measuring the dilatation of each compact edge, the cone
$C(R):=V(R)\cap \R^\dn_{\geq 0}$ parametrizes exactly the Minkowski summands
of positive multiples of $Q(R)$.
Hence, we will call elements of $V(R)$ ``generalized Minkowski summands'';
they may have edges of negative length.
(See \cite{versal}, Lemma (2.2) for a discussion of the compact case.)
The vector space $W(R)$ provides
coordinates  $s_\aj $ for each vertex $\bar{a}^\aj \in Q(R)\setminus N$, i.e.\
$\langle a^\aj ,R\rangle \geq 2$.\\
\par

%%%%%%%%%%
% (T1.7)
%%%%%%%%%

\neu{T1-7}
To each compact edge $d^{\aj\ak}=\overline{\bar{a}^\aj \bar{a}^\ak }$ we
assign a set of equations $\ksp_{\aj\ak}$ which act on elements of
$V(R)\oplus W(R)$. These sets are of one of the following three types:
\begin{itemize}
\item[(0)]
$\ksp_{\aj\ak}=\emptyset$,
\item[(1)]
$\ksp_{\aj\ak} = \{ s_\aj -s_\ak=0\}$ provided both coordinates exist in $W(R)$,
set $\ksp_{\aj\ak}=\emptyset$ otherwise, or
\item[(2)]
$\ksp_{\aj\ak} = \{t_{\aj\ak}-s_\aj=0 ,\; t_{\aj\ak}-s_\ak=0\}$,
dropping equations that do not make sense.
\end{itemize}
Restricting $V(R)\oplus W(R)$ to the (at most) three coordinates
$t_{\aj\ak}$, $s_\aj$, $s_\ak$,
the actual choice of $\ksp_{\aj\ak}$ is made such that these equations yield a
subspace of dimension $1+\mbox{dim}\,T^1_{\langle a^\aj,a^\ak\rangle}(-R)$.
Notice that the dimension of $T^1(-R)$ for the two-dimensional quotient singularity
assigned to the plane cone $\langle a^\aj,a^\ak\rangle$ can be obtained
from Example \zitat{T1}{3}.\\
\par

{\bf Theorem:}
{\em
The infinitesimal deformations of $Y_\sigma$ in degree $-R$ equal
\[
T^1_Y(-R)=
\{ (\ku{t},\,\ku{s})\in V_{\C}(R)\oplus W_{\C}(R)\,|\;
(\ku{t},\,\ku{s}) \mbox{ fulfills the equations } \ksp_{\aj\ak}\}
\;\big/\; \C\cdot (\ku{1},\, \ku{1})\,.
\vspace{-1ex}
\]
}

In some sense, the vector space $V(R)$ (encoding Minkowski summands)
may be considered the main tool to describe infinitesimal deformations.
The elements of $W(R)$ can (depending on the type of the $\ksp_{\aj\ak}$'s)
be either additional parameters, or they provide conditions excluding
Minkowski summands not having some prescribed type.\\
If $Y$ is smooth in codimension two, then $\ksp_{\aj\ak}$ is always of type (2).
In particular, the variables $\ku{s}$ are completely determined by the
$\ku{t}$'s, and we obtain the\\
\par

{\bf Corollary:}
{\em If $Y$ is smooth in codimension two, then
$T^1_Y(-R)$ is contained in $V_{\C}(R) \big/ \,\C\cdot (\ku{1})$. It is built from those
$\ku{t}$ such that $t_{\aj\ak}=t_{\ak\al}$ whenever $d^{\aj\ak}$, $d^{\ak\al}$
are compact edges with a common non-lattice vertex $\bar{a}^\ak$ of $Q(R)$.\\
Thus, $T^1_Y(-R)$ equals the set of equivalence classes of those Minkowski
summands of $\R_{\geq 0}\cdot Q(R)$ that preserve up to homothety the stars
of non-lattice vertices of $Q(R)$.
}\\
\par

%%%%%%%%%%
% (T1.8)
%%%%%%%%%

\neu{T1-8}
{\bf Proof:}\quad (of previous theorem)\\
{\em Step 1:}\quad
From Theorem \zitat{T1}{2} we know that $T^1_Y(-R)$ equals the
complexification of the cohomology of the complex
\[
N_{\R}\rightarrow
\oplus_\aj  \left(\mbox{span}_{\R} E_\aj^R \right)^\ast
\rightarrow
\oplus_{\langle a^\aj ,a^\ak \rangle <\sigma}
\left(\mbox{span}_{\R} E^R_{\aj\ak} \right)^\ast\,.
\]
According to \zitat{T1}{2}, elements of
$\oplus_\aj  \left(\mbox{span}_{\R} E_\aj^R \right)^\ast$
can be represented by a family of
\[
b^\aj \in N_{\R}\; \mbox{ (if } \langle a^\aj ,R\rangle \geq 2)\quad
\mbox{ and }
\quad b^\aj \in N_{\R}\big/\R\cdot a^\aj\; \mbox{ (if }
\langle a^\aj ,R\rangle =1).
\]
Dividing by the image of $N_{\R}$ means to shift this family by common
vectors $b\in N_{\R}$.
On the other hand, the family $\{b^\aj \}$ has to map onto $0$ in the complex,
i.e.\ for each compact edge
$\overline{\bar{a}^\aj ,\bar{a}^\ak }<Q$ the functions $b^\aj $ and
$b^\ak $ are
equal on $\mbox{span}_{\R}E_{\aj\ak}^R$. Since
\[
(a^\aj ,a^\ak )^\bot \subseteq \mbox{span}_{\R}E_{\aj\ak}^R \subseteq
(\mbox{span}_{\R}E_\aj^R) \cap (\mbox{span}_{\R}E_\ak^R)\,,
\]
we immediately obtain the necessary condition
$b^\aj -b^\ak \in \R a^\aj  + \R a^\ak $.
However, the actual behavior of $\mbox{span}_{\R}E_{\aj\ak}^R$
will require a closer look (in the third step).\\
\par

{\em Step 2:}\quad
We introduce new ``coordinates'':
\begin{itemize}
\item
$\bar{b}^\aj := b^\aj -\langle b^\aj , R \rangle \,\bar{a}^\aj  \in R^\bot$,
being well defined even in the case $\langle a^\aj , R \rangle =1$;
\item
$s_\aj :=-\langle b^\aj , R\rangle$
for $\aj$ meeting $\langle a^\aj , R \rangle \geq 2$ (inducing an element
of $W(R)$).
\end{itemize}

The shift of the $b^\aj$ by an element $b\in N_{\R}$ (i.e.\
$(b^\aj )'=b^\aj +b$) appears in these new coordinates as
\[
\renewcommand{\arraystretch}{1.5}
\begin{array}{rcl}
(\bar{b}^\aj )' &=& (b^\aj )' - \langle (b^\aj )', R \rangle \,\bar{a}^\aj
\;=\; b^\aj  +b - \langle b^\aj ,R \rangle \,\bar{a}^\aj  -
\langle b,R \rangle\,\bar{a}^\aj
\vspace{-0.5ex}\\
&=& \bar{b}^\aj  + b-\langle b,R \rangle \,\bar{a}^\aj \,,\\
s_\aj '&=& -\langle (b^\aj )',R \rangle
\;=\; s_\aj -\langle b,R\rangle\,.
\end{array}
\]
In particular, an element $b\in R^\bot$ does not change the $s_\aj$,
but shifts the points $\bar{b}^\aj$ inside the hyperplane $R^\bot$. Hence,
the set of the $\bar{b}^\aj $ should be considered modulo translation
inside $R^\bot$ only.\\
On the other hand, the condition $b^\aj -b^\ak \in \R a^\aj  + \R a^\ak $
changes into
$\bar{b}^\aj -\bar{b}^\ak \in \R \bar{a}^\aj  + \R \bar{a}^\ak $ or even
$\bar{b}^\aj -\bar{b}^\ak \in \R (\bar{a}^\aj  - \bar{a}^\ak )$ (consider the values
of $R$). Hence, the $\bar{b}^\aj $'s form the vertices of an at least
generalized Minkowski summand of $Q(R)$. Modulo translation, this summand
is completely described by the dilatation factors $t_{\aj\ak}$ obtained from
\[
\bar{b}^\aj -\bar{b}^\ak  = t_{\aj\ak}\cdot (\bar{a}^\aj  - \bar{a}^\ak )\,.
\]
Now, the remaining part of
the action of $b\in N_{\R}$ comes down to an action of
$\langle b,R\rangle\in\R$ only:
\[
\begin{array}{rcl}
t_{\aj\ak}' &=& t_{\aj\ak} - \langle b,R \rangle \quad \mbox{ and}\\
s_\aj ' &=& s_\aj  - \langle b,R \rangle \,, \mbox{ as we already know}.
\end{array}
%\vspace{2ex}
\]
Up to now, we have found that
$T^1_Y(-R)\subseteq V_{\C}(R)\oplus W_{\C}(R)/(\ku{1},\ku{1})$.\\
\par

{\em Step 3:}\quad
Actually, the elements $b^\aj $ and $b^\ak $ have to be equal on
$\mbox{span}_{\R} E^R_{\aj\ak}$, which may be a larger space than just
$(a^\aj ,a^\ak )^\bot$.
To measure the difference we consider the factor
$\mbox{span}_{\R} E^R_{\aj\ak}\big/ (a^\aj ,a^\ak )^\bot$
contained in the two-dimensional vector space
$M_{\R}\big/ (a^\aj ,a^\ak )^\bot =\mbox{span}_{\R}(a^\aj ,a^\ak )^\ast$.
Since this factor coincides with the set $\mbox{span}_{\R}E^{\bar{R}}_{\aj\ak}$
assigned to the two-dimensional cone
$\langle a^\aj ,a^\ak \rangle \subseteq \mbox{span}_{\R}(a^\aj ,a^\ak )$,
where $\bar{R}$
denotes the image of $R$ in $\mbox{span}_{\R}(a^\aj ,a^\ak )^\ast$,
we may assume
that $\sigma=\langle a^1, a^2\rangle$ (i.e.\ $\aj=1,\,\ak=2$)
represents a two-dimensional cyclic
quotient singularity. In particular, we only need to discuss the three cases
(i)-(iii) from Example \zitat{T1}{3}:\\
In (i) and (ii) we have $\mbox{span}_{\R} E^R_{12}=0$, i.e.\ no additional
equation is needed. This means $\ksp_{12}=\emptyset$
is of type (0). On the other hand, if $T^1_Y=0$, then the
vector space $\R^3_{(t_{12},s_1,s_2)}\big/\R\cdot (\ku{1})$ has to be killed
by identifying the three variables $t_{12}$, $s_1$, and $s_2$;
we obtain type (2).\\
Case (iii) provides $\mbox{span}_{\R} E^R_{12}=\R\cdot R$. Hence, as an
additional condition we obtain that $b^1$ and $b^2$ have to be equal on $R$.
By the definition of $s_\aj$ this means $s_1=s_2$, and $\ksp_{12}$ has
to be of type (1).
\hfill$\Box$\\
\par

%%%%%%%%%%%%
%
%  Genuine Deformations
%
%%%%%%%%%%%%%
\section{Genuine deformations}\label{Gd}

%%%%%%%%%%
% (Gd.1)
%%%%%%%%

\neu{Gd-1}
In \cite{Tohoku} we have studied so-called toric deformations in a given
multidegree $-R\in M$. They are genuine deformations in the sense that they are
defined over smooth parameter spaces;
they are characterized by the fact that the total spaces
together with the embedding of the special fiber still belong to the toric
category. Despite the fact they look so special, it seems that toric deformations
cover a big part of the versal deformation of $Y_\sigma$. They do only
exist in negative degrees (i.e.\ $R\in\sigma^{\veee}\cap M$), but here they
form a kind of skeleton. If $Y_\sigma$ is an isolated toric Gorenstein
singularity, then toric deformations even provide all irreducible components
of the versal deformation (cf.\ \cite{versal}).\\
After a quick reminder of the idea of this construction, we
describe the Kodaira-Spencer map of toric deformations in terms of the new
$T^1_Y$-formula presented in \zitat{T1}{2}. It is followed by the investigation
of non-negative degrees: If $R\notin\sigma^{\veee}\cap M$, then we are still
able to construct genuine deformations of $Y_\sigma$; but they are no longer toric.\\
\par

%%%%%%%%%%
% (Gd.2)
%%%%%%%%

\neu{Gd-2}
Let $R\in\sigma^{\veee}\cap M$. Then, following \cite{Tohoku} \S 3,
toric $\sm$-parameter deformations of $Y_\sigma$ in degree $-R$ correspond
to splittings of $Q(R)$ into a Minkowski sum
\vspace{-0.5ex}
\[
Q(R) \,=\, Q_0 + Q_1 + \dots +Q_\sm
\vspace{-1ex}
\]
meeting the following conditions:
\begin{itemize}
\item[(i)]
$Q_0\subseteq \A{R}{}$ and $Q_1,\dots,Q_\sm\in \A{R}{0}$ are polyhedra with $Q(R)^\infty$
as their common cone of unbounded directions.
\item[(ii)]
Each supporting hyperplane $t$ of $Q(R)$ 
% may be given by some element of $(Q(R)^\infty)^{\veee}\subseteq \A{R}{0}^\ast$ and
defines faces
$F(Q_0,t),\dots, F(Q_\sm,t)$ of the indicated polyhedra; their Minkowski sum
equals $F\big(Q(R),t\big)$.
With at most one exception (depending on $t$), these faces should contain
lattice vertices, i.e.\ vertices belonging to $N$.
\end{itemize}

{\bf Remark:}
In \cite{Tohoku} we have distinguished between the case of primitive
and non-primitive
elements $R\in M$: If $R$ is a multiple of some element of $M$, then $\A{R}{}$
does not
contain lattice points at all. In particular, condition (ii) just means that
$Q_1,\dots,Q_\sm$ have to be lattice polyhedra.\\
On the other hand, for primitive $R$, the $(\sm+1)$ summands $Q_\si$ have
equal rights
and may be put into the same space $\A{R}{}$. Then, their Minkowski sum has to
be interpreted inside this affine space.\\
\par

If a Minkowski decomposition is given, {\em how do we obtain the assigned toric
deformation?}\\
Defining $\tN:= N\oplus \Z^\sm$ (and $\tM:=M\oplus \Z^\sm$),
we have to embed the summands as $(Q_0,\,0)$,
$(Q_1,\,e^1),\dots, (Q_\sm,\,e^\sm)$ into the vector space
$\tN_{\R}$; $\{e^1,\dots,e^\sm\}$ denotes
the canonical basis of $\Z^\sm$. Together with $(Q(R)^\infty,\,0)$, these
polyhedra generate
a cone $\tsigma\subseteq \tN$ containing $\sigma$ via
$N\hookrightarrow \tN$, $a\mapsto (a;\langle a,R\rangle,\dots,\langle a,R\rangle)$.
Actually, $\sigma$ equals $\tsigma\cap N_{\R}$, and we obtain an inclusion
$Y_\sigma\hookrightarrow X_{\tsigma}$ between the associated toric varieties.\\
On the other hand, $[R,0]:\tN\to \Z$ and $\mbox{pr}_{\Z^\sm}:\tN\to\Z^\sm$ induce
regular
functions $f:X_{\tsigma}\to \C$ and $(f^1,\dots,f^\sm):X_{\tsigma}\to \C^\sm$,
respectively. The resulting map $(f^1-f,\dots,f^\sm-f):X_{\tsigma}\to \C^\sm$ is flat
and has $Y_\sigma\hookrightarrow X_{\tsigma}$ as special fiber.\\
\par

%%%%%%%%%%
% (Gd.3)
%%%%%%%%

\neu{Gd-3}
Let $R\in\sigma^{\veee}\cap M$ and $Q(R) = Q_0 + \dots +Q_\sm$ be a decomposition
satisfying (i) and (ii) mentioned above. Denote by $(\bar{a}^\aj)_\si$ the vertex
of $Q_i$ induced from $\bar{a}^\aj\in Q(R)$, i.e.\
$\bar{a}^\aj=(\bar{a}^\aj)_0 + \dots + (\bar{a}^\aj)_\sm$.\\
\par

{\bf Theorem:}
{\em
The Kodaira-Spencer map of the corresponding toric deformation
$X_{\tsigma}\to \C^\sm$ is
\[
\varrho: \C^\sm\,=\, T_{\C^\sm,0}  \longrightarrow T^1_Y(-R)
\subseteq V_{\C}(R)\oplus W_{\C}(R)\big/\C\cdot (\ku{1},\ku{1})
\]
sending $e^\si $ onto the pair $[Q_\si ,\,\ku{s}^\si ]\in V(R)\oplus W(R)$
($\si=1,\dots,\sm$) with
\[
s^\si_\aj :=\left\{ \begin{array}{ll}
0 & \mbox{ if the vertex } (\bar{a}^\aj)_\si  \mbox{ of } Q_\si
\mbox{ belongs to the lattice } N\\
1 & \mbox{ if } (\bar{a}^\aj)_\si  \mbox{ is not a lattice point.}
\end{array}\right.
\]
}
\par

{\bf Remark:}
Setting $e^0:=-(e^1+\dots +e^\sm )$, we obtain
$\varrho (e^0)= [Q_0,\, \ku{s}^0]$ with $\ku{s}^0$ defined similar to $\ku{s}^\si $
in the previous theorem.\\
\par

%%%%%%%%%%
% (Gd.4)
%%%%%%%%

\neu{Gd-4}
{\bf Proof} (of previous theorem):
We would like to derive the above formula for the Kodaira-Spencer map from the
more technical one presented in \cite{Tohoku}, Theorem (5.3).
% (For the $V(R)$-part, this has more or less been done already in (5.5)
% of the same paper.)
Under additional use of \cite{T2} (6.1),
the latter one describes $\varrho(e^\si)\in
T^1_Y(-R)=H^1\big(\mbox{span}_{\C}(E^R)_\bullet^\ast\big)$
in the following way:\\
Let $E=\{r^1,\dots,r^w\}\subseteq \sigma^{\veee}\cap M$. Its elements
may be lifted via $\tM\surj M$
% $(r,\ku{\eta})\mapsto r +(\eta_1+\dots +\eta_\sm)\,R$
to $\tilde{r}^v\in\tsigma^{\veee}\cap\tM$ ($v=1,\dots,w$); denote their $\si$-th
entry of the $\Z^\sm$-part by $\eta^v_\si$, respectively.
Then, given elements $v^\aj\in \mbox{span} E_\aj^R$, we may represent them
as $v^\aj=\sum_v q^j_v\,r^v$ ($q^\aj\in\Z^{E_\aj^R}$), and $\varrho(e^\si)$
assigns to $v^\aj$ the integer $-\sum_v q^j_v\,\eta_\si^v$.
Using our notation from \zitat{T1}{8} for $\varrho(e^\si)$, this means that
$b^\aj$ sends elements $r^v\in E_\aj^R$ onto $-\eta_\si^v\in\Z$. \\
By construction of $\tsigma$, we have inequalities
\[
\Big\langle \big( (\bar{a}^\aj)_0,\,0\big),\, \tilde{r}^v \Big\rangle \geq 0
\quad\mbox{ and }\quad
\Big\langle \big( (\bar{a}^\aj)_\si,\,e^i\big),\, \tilde{r}^v \Big\rangle \geq 0
\quad (\si=1,\dots,\sm)
\]
summing up to $\big\langle \bar{a}^\aj,\, r^v \big\rangle =
\big\langle \big( \bar{a}^\aj,\,\ku{1}\big),\, \tilde{r}^v \big\rangle \geq 0$.
On the other hand, the fact $r^v\in E_\aj^R$ is equivalent to
$\big\langle \bar{a}^\aj,\, r^v \big\rangle <1$. Hence, whenever
$(\bar{a}^\aj)_\si\in Q_\si$ belongs to the lattice, the
corresponding inequality ($\si=0,\dots,\sm$) becomes an equality.
With at most one exception, this always has to be the case. Hence,
\[
\Big\langle (\bar{a}^\aj)_\si,\, r^v \Big\rangle + \eta_\si^v
= \; \left\{
\begin{array}{ll}
0& \mbox{ if } (\bar{a}^\aj)_\si \in N\\
\langle \bar{a}^\aj,\, r^v \rangle & \mbox{ if } (\bar{a}^\aj)_\si \notin N
\end{array}
\right.
\quad(\si=1,\dots,\sm)
\]
meaning that $b^\aj= (\bar{a}^\aj)_\si\,$ or
$\,b^\aj= (\bar{a}^\aj)_\si - \bar{a}^\aj$, respectively.
By the
definitions of $\bar{b}^\aj$ and $s_\aj$ given in \zitat{T1}{8}, we are done.
\hfill$\Box$\\
\par

%%%%%%%%%%
% (Gd.5)
%%%%%%%%

\neu{Gd-5}
Now we treat the case of non-negative degrees; let $R\in M\setminus \sigma^{\veee}$.
The easiest way to solve a problem is to change the question until there is no problem
left. We can do so by changing our cone $\sigma$ into some $\tau^R$ such that the
degree $-R$ becomes negative. We define
\[
\tau:=\tau^R:= \sigma \cap \,[R\geq 0]\quad
\mbox{ that is } \quad
\tau^{\veee}=\sigma^{\veee}+\R_{\geq 0}\cdot R\,.
\]
The cone $\tau$ defines an affine toric variety $Y_\tau$. Since
$\tau\subseteq\sigma$, it comes with a map $g:Y_\tau\to Y_\sigma$, i.e.\
$Y_\tau$ is an open part of a modification of $Y_\sigma$. The important
observation is
\[
\renewcommand{\arraystretch}{1.5}
\begin{array}{rcccl}
\tau \cap \,[R=0] &=& \sigma\cap\, [R=0] &=& Q(R)^\infty\quad \mbox{ and}\\
\tau \cap \,[R=1] &=& \sigma\cap\, [R=1] &=& Q(R)\;,
\end{array}
\]
implying $T^1_{Y_\tau}(-R)=T^1_{Y_\sigma}(-R)$ by Theorem \zitat{T1}{7}. Moreover,
even the genuine toric deformations $X_\ttau\to\C^\sm$ of $Y_\tau$ carry over to
$\sm$-parameter (non-toric) deformations $X\to\C^\sm$ of $Y_\sigma$:\\
\par

{\bf Theorem:}
{\em
Each Minkowski decomposition $Q(R) = Q_0 + Q_1 + \dots +Q_\sm$
satisfying (i) and (ii) of \zitat{Gd}{2}
provides an $\sm$-parameter deformation $X\to\C^\sm$ of $Y_\sigma$. Via some
birational map $\tilde{g}:X_\ttau\to X$ it is compatible with the
toric deformation $X_\ttau\to \C^\sm$ of $Y_\tau$ presented in
\zitat{Gd}{2}.
% They are arranged in the following commutative diagram:
\[
\dgARROWLENGTH=0.4em
\begin{diagram}
\node[2]{\C^\sm}
\arrow[3]{e,t}{\mbox{\kf id}}
\node[3]{\C^\sm}\\
\node{X_\ttau}
\arrow{ne}
\arrow[3]{e,t}{\tilde{g}}
\node[3]{X}
\arrow[2]{e}
\arrow{ne}
\node[2]{Z_{\tsigma}}\\[2]
\node{Y_\tau}
\arrow[2]{n}
\arrow[3]{e,t}{g}
\node[3]{Y_\sigma}
\arrow[2]{n}
\arrow[2]{ne}
\end{diagram}
\]
The total space $X$ is not toric anymore, but it sits via birational maps between 
$X_\ttau$ and some affine toric variety $Z_{\tsigma}$
still containing $Y_\sigma$ as a closed subset. 
}\\
\par

%%%%%%%%%%
% (Gd.6)
%%%%%%%%

\neu{Gd-6}
{\bf Proof:}
First, we construct $\tN$, $\tM$, and $\ttau\subseteq\tN_{\R}$ by the recipe
of \zitat{Gd}{2}. In particular, $N$ is contained in $\tN$, and the projection
$\pi:\tM\to M$ sends $[r;g_1,\dots,g_\sm]$ onto $r+(\sum_\si g_\si)\,R$.
Defining $\tsigma:= \ttau + \sigma$ 
(hence $\tsigma^{\veee}=\ttau^{\veee}\cap \pi^{-1}(\sigma^{\veee})$), we obtain
the commutative diagram
\[
\dgARROWLENGTH=0.5em
\begin{diagram}
\node{\C[\ttau^{\veee}\cap \tM]}
\arrow{s,r}{\pi}
\node[3]{\C[\tsigma^{\veee}\cap \tM]}
\arrow{s,r}{\pi}
\arrow[3]{w}\\
\node{\C[\tau^{\veee}\cap M]}
\node[3]{\C[\sigma^{\veee}\cap M]}
\arrow[3]{w}
\end{diagram}
\]
with surjective vertical maps. The canonical elements
$e_1,\dots,e_\sm\in\Z^\sm\subseteq\tM$ together with $[R;0]\in\tM$ are preimages
of $R\in M$. Hence, the corresponding monomials
$x^{e_1},\dots,x^{e_\sm},x^{[R,0]}$ in the semigroup algebra
$\C[\ttau^{\veee}\cap \tM]$ (called $f^1,\dots,f^\sm,f$ in \zitat{Gd}{2})
map onto $x^R\in\C[\tau^{\veee}\cap M]$ which is not regular on $Y_\sigma$.
We define $Z_\tsigma$ as the affine toric variety assigned to $\tsigma$ and
$X$ as
\[
X:=\mbox{Spec}\,B \quad \mbox{ with } \quad
B:=\C[\tsigma^{\veee}\cap\tM][f^1-f,\dots,f^\sm-f]\subseteq
\C[\ttau^{\veee}\cap\tM]\,.
\]
That means, $X$ arises from $X_\ttau$ by eliminating all variables except
those lifted from $Y_\sigma$ or the deformation parameters themselves.
By construction of $B$, the vertical algebra homomorphisms $\pi$ induce
a surjection $B\surj \C[\sigma^{\veee}\cap M]$.\\
\par

{\em Lemma:} Elements of $\C[\ttau^{\veee}\cap\tM]$ may uniquely be written
as sums
\vspace{-1ex}
\[
\sum_{(v_1,\dots,v_\sm)\in\N^\sm} c_{v_1,\dots,v_\sm}\cdot
(f^1-f)^{v_1}\cdot\dots\cdot (f^\sm-f)^{v_\sm}
\vspace{-1.5ex}
\]
with $c_{v_1,\dots,v_\sm}\in\C[\ttau^{\veee}\cap\tM]$
such that $s-e_\si\notin \ttau^{\veee}$ ($\si=1,\dots,\sm$) for any
of its monomial terms $x^s$. Moreover, those sums belong to the subalgebra $B$,
if and only if their coefficients $c_{v_1,\dots,v_\sm}$ do.
\vspace{1ex}\\
{\em Proof:
(a) Existence.}
Let $s-e_\si\in\ttau^{\veee}$ for some $s,\si$. Then, with
% for some term in some coefficient $c_{v_1,\dots,v_\sm}$
$s^\prime:=s-e_\si+[R,0]$ we obtain
\vspace{-1ex}
\[
x^s = x^{s^\prime} + x^{s-e_\si}\,(x^{e_\si}-x^{[R,0]}) =
      x^{s^\prime} + x^{s-e_\si}\,(f^\si-f)\,.
\vspace{-1ex}
\]
Since $e_\si=1$ and $[R,0]=0$ if evaluated on $(Q_\si,e^\si)\subseteq\ttau$,
this process eventually stops.
\vspace{1ex}\\
{\em (b) $B$-Membership.}
For the previous reduction step we have to show that if
$s\in\C[\tsigma^{\veee}\cap\tM]$, then the same is true for $s^\prime$ and
$s-e_\si$.
Since $\pi(s^\prime)=\pi(s)\in\sigma^{\veee}$, this is clear for $s^\prime$.
It remains to
check that $\pi(s-e_\si)\in\sigma^{\veee}$.
Let $a\in\sigma$ be an arbitrary test element; we distinguish two cases:\\
Case 1: $\langle a,R\rangle\geq 0$. Then $a$ belongs to the subcone
$\tau$, and $\pi(s-e_\si)\in \tau^{\veee}$ yields
$\langle a, \pi(s-e_\si)\rangle\geq 0$.\\
Case 2: $\langle a,R\rangle\leq 0$. This fact implies
$\langle a, \pi(s-e_\si)\rangle = \langle a, s\rangle - \langle a,R\rangle
\geq \langle a, s\rangle \geq 0$.
\vspace{1ex}\\
{\em (c) Uniqueness.} Let $p:=\sum c_{v_1,\dots,v_\sm}\cdot
(f^1-f)^{v_1}\cdot\dots\cdot (f^\sm-f)^{v_\sm}$ (meeting the above
conditions) be equal to $0$ in $\C[\ttau^{\veee}\cap\tM]$. Using the projection
$\pi:\tM\to M$, everything becomes $M$-graded. Since the factors
$(f^\si-f)$ are homogeneous (of degree $R$), we may assume this fact also
for $p$, hence for its coefficients $c_{v_1,\dots,v_\sm}$.
\vspace{0.5ex}\\
{\em Claim:} These coefficients are just monomials. Indeed, if
$s,s^\prime\in\ttau^{\veee}$ had the same image via $\pi$, then we
could assume that some $e_\si$-coordinate of $s^\prime$
would be smaller than that
of $s$. Hence, $s-e_\si$ would still be equal to $s$ on $(Q_0,0)$ and on
any $(Q_\sj,e^\sj)$ ($\sj\neq\si$), but even greater or equal than $s^\prime$
on $(Q_\si,e^\si)$. This would imply $s-e_\si\in\ttau^{\veee}$,
contradicting our assumption for $p$.
\vspace{0.5ex}\\
Say $c_{v_1,\dots,v_\sm}=\lambda_{v_1,\dots,v_\sm}\,x^\bullet$;
we use the projection $\tM\to\Z^\sm$ for carrying $p$ into the ring
$\C[\Z^\sm]=\C[y_1^{\pm 1},\dots,y_\sm^{\pm 1}]$. The elements
$x^\bullet$, $f^\si$, $f$ map onto $y^\bullet$, $y_\si$, and $1$,
respectively. Hence, $p$ turns into
\vspace{-1ex}
\[
\bar{p}=\sum_{(v_1,\dots,v_\sm)\in\N^\sm} \lambda_{v_1,\dots,v_\sm}\cdot
y^\bullet\cdot(y_1-1)^{v_1}\cdot\dots\cdot (y_\sm-1)^{v_\sm}\,.
\vspace{-1ex}
\]
By induction through $\N^\sm$, we obtain that vanishing of $\bar{p}$
implies the vanishing of its coefficients: Replace $y_\si-1$ by $z_\si$,
and take partial derivatives.
\hfill$(\Box)$\\
\par

Now, we can easily see that $X\to\C^\sm$ is flat and has $Y_\sigma$ as
special fiber:
The previous lemma means that for $\sk=0,\dots,\sm$ we have inclusions
\[
{\kd B}\big/_{\kd (f^1-f,\dots,f^\sk-f)}
\quad\raisebox{-0.5ex}{$\hookrightarrow\;$}\quad
{\kd \C[\ttau^{\veee}\cap\tM]\,}\big/_{\kd (f^1-f,\dots,f^\sk-f)}\,.
\]
The values $\sk < \sm$ yield that $(f^1-f,\dots,f^\sm-f)$ forms a regular
sequence even in the subring $B$, meaning that $X\to\C^\sm$ is flat.
With $\sk=\sm$ we obtain that the surjective map
$B/(f^1-f,\dots,f^\sm-f)\to\C[\sigma^{\veee}\cap M]$ is also injective.
\hfill$\Box$\\
\par

%%%%%%%%%%%%
%
%  Gorenstein Singularities
%
%%%%%%%%%%%%%
\section{Three-dimensional toric Gorenstein singularities}\label{3G}

%%%%%%%%%%
% (3G.1)
%%%%%%%%

\neu{3G-1}
By \cite{Ish}, Theorem (7.7),
toric Gorenstein singularities always arise from the following construction:
Assume we are given a {\em lattice polytope $P\subseteq \R^\pn$}.
We embed the whole space
(including $P$) into height one of $N_{\R}:=\R^\pn\oplus\R$ and take for
$\sigma$ the
cone generated by $P$; denote by $M_{\R}:=(\R^\pn)^\ast\oplus\R$
the dual space and by $N$, $M$ the natural lattices.
Our polytope $P$ may be recovered from $\sigma$ as
\[
P\, =\, Q(R^\ast)\subseteq\A{R^\ast}{}\quad
\mbox{ with} \quad R^\ast:=[\ku{0},1]\in M\,.
\hspace{-3em}
\raisebox{-35mm}{
\unitlength=0.5mm
\linethickness{0.6pt}
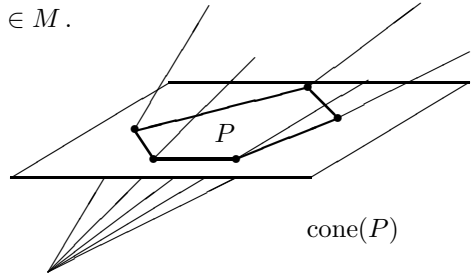
\begin{picture}(130,80)
\thinlines
\put(0.00,30.00){\line(1,0){80.00}}
\put(0.00,30.00){\line(5,3){42.00}}
\put(42.00,55.33){\line(1,0){80.00}}
\put(122.00,55.33){\line(-5,-3){42.00}}
\put(10.00,5.00){\line(3,5){15.00}}
\put(33.00,43.00){\line(3,5){19.67}}
\put(10.00,5.00){\line(1,1){25.00}}
\put(38.00,35.00){\line(1,1){27.00}}
\put(10.00,5.00){\line(5,3){42.00}}
\put(60.00,35.00){\line(5,3){35.00}}
\put(10.00,5.00){\line(2,1){50.00}}
\put(87.00,46.00){\line(2,1){35.00}}
\put(10.00,5.00){\line(4,3){33.00}}
\put(79.00,54.00){\line(4,3){30.00}}
\put(57.00,42.00){\makebox(0,0)[cc]{$P$}}
\put(91.00,16.00){\makebox(0,0)[cc]{$\mbox{cone}(P)$}}
\put(33.00,43.00){\circle*{2.00}}
\put(38.00,35.00){\circle*{2.00}}
\put(60.00,35.00){\circle*{2.00}}
\put(87.00,46.00){\circle*{2.00}}
\put(79.00,54.00){\circle*{2.00}}
\thicklines
\put(32.67,43.00){\line(2,-3){5.33}}
\put(38.00,35.00){\line(1,0){22.00}}
\put(60.00,35.00){\line(5,2){27.00}}
\put(87.00,46.00){\line(-1,1){8.00}}
\put(79.00,54.00){\line(-4,-1){46.00}}
\end{picture}
}
\]
The fundamental generators $a^1,\dots,a^\an\in\G{R^\ast}{}$ of $\sigma$
coincide with the vertices of $P$. (This involves a slight abuse of notation;
we use the same symbol $a^\aj$ for both $a^\aj\in\Z^\pn$ and $(a^\aj,1)\in M$.)\\
If $\ko{a^\aj a^\ak}$ forms an edge of the polytope,
we denote by $\ell(\aj,\ak)\in\Z$ its ``length'' induced from the
lattice structure $\Z^\pn\subseteq\R^\pn$. Every edge provides a
two-codimensional singularity of $Y_\sigma$ with transversal type
A$_{\ell(\aj,\ak)-1}$. In particular, $Y_\sigma$ is smooth in codimension
two if and only if all edges of $P$ are primitive, i.e.\ have length
$\ell=1$.\\
\par

%%%%%%%%%%
% (3G.2)
%%%%%%%%

\neu{3G-2}
As usual, we fix some element $R\in M$. From \zitat{T1}{6} we know what the
vector spaces $V(R)$ and $W(R)$ are; we introduce the subspace
\[
V^\prime(R):=\{\ku{t}\in V(R)\,|\; t_{\aj\ak}\neq 0 \mbox{ implies }
1\leq \langle a^\aj,R \rangle = \langle a^\ak,R \rangle \leq \ell(\aj,\ak)\}
\]
representing Minkowski summands of $Q(R)$ that have killed any compact edge
{\em not} meeting the condition
$\langle a^\aj,R \rangle = \langle a^\ak,R \rangle \leq \ell(\aj,\ak)$.\\
\par

{\bf Theorem:}
{\em
For $T^1_Y(-R)$, there are two different types of $R\in M$ to distinguish:
\begin{itemize}
\item[(i)]
If $R\leq 1$ on $P$ (or equivalently $\langle a^\aj,R\rangle\leq 1$ for
$\aj=1,\dots,\an$), then $T^1_Y(-R)=V_{\C}(R)\big/(\ku{1})$. Moreover,
concerning Minkowski summands, we may replace the polyhedron
$Q(R)$ by its compact part $P\cap [R=1]$ (being a face of $P$).
\item[(ii)]
If $R$ does not satisfy the previous condition, then
$T^1_Y(-R)=V^\prime(R)$.
\vspace{1ex}
\end{itemize}
}

{\bf Proof:}
The first case follows from Theorem \zitat{T1}{7} just because $W(R)=0$.
For (ii), let us assume there are vertices $a^\aj$ contained in the affine
half space $[R\geq 2]$.
They are mutually connected inside this half space via paths along edges of $P$.\\
The two-dimensional cyclic quotient singularities corresponding to edges
$\ko{a^\aj a^\ak}$ of $P$ are Gorenstein themselves. In the language of
Example \zitat{T1}{3} this means $w=2$, and we obtain
\[
\mbox{dim}\, T^1_{\langle a^\aj,a^\ak \rangle} (-R)\,=\,
\left\{ \begin{array}{ll}
1 & \mbox{ if } \langle a^\aj,R \rangle = \langle a^\ak,R \rangle
=2,\dots,\ell(\aj,\ak)
\quad\mbox{(case (iii) in \zitat{T1}{3})}\\
0 & \mbox{ otherwise.}
\end{array}\right.
\]
In particular, $T^1_{\langle a^\aj,a^\ak \rangle} (-R)$ cannot be two-dimensional,
and (using the notation of \zitat{T1}{7})
the equations $s_\aj -s_\ak=0$ belong to $\ksp_{\aj\ak}$ whenever
$\langle a^\aj,R \rangle ,\, \langle a^\ak,R \rangle\geq 2$.
This means for elements of
\[
T^1_Y\subseteq  \Big(V_{\C}(R)\oplus W_{\C}(R)\Big)\Big/ \C\cdot (\ku{1},\ku{1})
\]
that all entries
of the $W_{\C}(R)$-part have to be mutually equal, or even zero after dividing by
$\C\cdot (\ku{1},\ku{1})$.
Moreover, if not both $\langle a^\aj,R \rangle$ and $\langle a^\ak,R \rangle$
equal one, vanishing of $T^1_{\langle a^\aj,a^\ak \rangle} (-R)$
implies that $\ksp_{\aj\ak}$ also contains the equation $t_{\aj\ak} -s_\bullet=0$.
\hfill$\Box$\\
\par

{\bf Corollary:}
{\em
Condition \zitat{Gd}{2}(ii) to build genuine deformations becomes easier for
toric Gorenstein singularities: $Q_1,\dots,Q_\sm$ just have to be lattice
polyhedra.
\vspace{-1ex}}\\
\par

{\bf Proof:}
If $R\leq 1$ on $P$, then $Q(R)$ itself is a lattice polyhedron. Hence,
condition (ii) automatically comes down to this simpler form.\\
In the second case, there is some $W(R)$-part involved in $T^1_Y(-R)$.
On the one hand, it
indicates via the Kodaira-Spencer map which
vertices of which polyhedron $Q_\si$ belong to the lattice. On the other,
we have observed in the previous proof that the entries of $W(R)$
are mutually equal. This implies exactly our claim.
\hfill$\Box$\\
\par

%%%%%%%
% (3G.3)
%%%%%%%%

\neu{3G-3}
In accordance with the title of the section, we focus now on {\em plane lattice polygons
$P\subseteq \R^2$}. The vertices $a^1,\dots,a^\an$ are arranged in a cycle.
We denote by $d^\aj:=a^{\aj+1}-a^\aj\in \G{R^\ast}{0}$ the edge going from $a^\aj$
to $a^{\aj+1}$, and by $\ell(\aj):=\ell(\aj,\aj+1)$ its length ($\aj\in\Z/\an\Z$).\\
Let $s^1,\dots,s^\an$ be the fundamental generators of the dual cone $\sigma^{\veee}$
such that $\sigma\cap(s^\aj)^\bot$ equals the face spanned by
$a^\aj, a^{\aj+1}\in\sigma$. In particular, skipping the last coordinate of
$s^\aj$ yields the (primitive) inner normal vector at the edge $d^\aj$ of $P$.
\vspace{-1ex}\\
\par

{\bf Remark:}
Just for convenience of those who prefer living in $M$ instead of $N$, we
show how to see the integers $\ell(\aj)$ in the dual world:
Choose a fundamental generator $s^\aj$ and denote by $\kr, \kr^\prime\in M$
the closest (to $s^\aj$) elements from the Hilbert bases of the two adjacent
faces of $\sigma^{\veee}$, respectively. Then, $\{R^\ast,s^\aj\}$
together with either $\kr$ or $\kr^\prime$ form a basis of the lattice $M$, and
$(\kr+\kr^\prime)-\ell(\aj)\,R^\ast$ is a positive multiple of $s^\aj$.
See the figure in \zitat{3G}{7}.\\
\par

In the very special case of plane lattice polygons (or three-dimensional toric
Gorenstein singularities), we can describe $T^1_Y$ and the genuine deformations (for
fixed $R\in M$) explicitly. First, we can easily spot the degrees carrying
infinitesimal deformations:
\vspace{-1ex}\\
\par

{\bf Theorem:}
{\em
In general (see the upcoming exceptions), $T^1_Y(-R)$ is non-trivial only for
\begin{itemize}
\item[(1)]
$R=R^\ast\,$ with $\,\mbox{\em dim}\,T^1_Y(-R)=\an-3$;
\item[(2)]
$R= qR^\ast$ ($q\geq 2)\,$ with
$\,\mbox{\em dim}\,T^1_Y(-R)= \mbox{\em max}\,\{0\,;\;
\#\{\aj\,|\; q\leq \ell(\aj)\}-2\,\}$, and
\item[(3)]
$R=qR^\ast - p\,s^\aj\,$ with $\,2\leq q\leq \ell(\aj)$ and
$p\in\Z$ sufficiently large such that $R\notin\mbox{\rm int}(\sigma^{\veee})$.
In this case, $T^1_Y(-R)$ is one-dimensional.
\end{itemize}
Additional degrees
% yielding a non-trivial $T^1_Y(-R)$
\vspace{-1ex}
exist only in the following two (overlapping) exceptional cases:
\begin{itemize}
\item[(4)]
Assume $P$ contains a pair of parallel edges $d^\aj$, $d^\ak$, both longer
than every other edge. Then $\mbox{\rm dim}\,T^1_Y(-q\,R^\ast)=1$ for
$\mbox{\rm max}\{\ell(\al)\,|\;\al\neq\aj,\ak\}<q\leq
\mbox{\rm min}\{\ell(\aj),\ell(\ak)\}$.
\vspace{-1ex}
\item[(5)]
Assume $P$ contains a pair of parallel edges $d^\aj$, $d^\ak$ with distance
$d$ ($d:=\langle a^\aj, s^\ak\rangle = \langle a^\ak,s^\aj\rangle$).
If $\ell(\ak)>d \;(\geq \mbox{\rm max}\{\ell(\al)\,|\;\al\neq\aj,\ak\})$, then
$\mbox{\rm dim}\,T^1_Y(-R)=1$ for
$R=qR^\ast +p\,s^\aj$ with
$1\leq q\leq\ell(\aj)$ and $1\leq p\leq \big(\ell(\ak)-q\big)/d$.
\end{itemize}
}
\par

The cases (1), (2), (4), and (5) yield at most finitely many
(negative) $T^1_Y$-degrees. Type (3) consists of $\ell(\aj)\!-\!1$ infinite series
to any vertex $a^\aj\in P$, respectively;
up to maybe the leading elements ($R$ might sit on
$\partial\sigma^{\veee}$), they contain only non-negative degrees.\\
\par

{\bf Proof:}
The previous claims are straight consequences of Theorem \zitat{3G}{2}.
Hence, the following short remark should be sufficient: The condition
$\langle a^\aj, R\rangle =\langle a^{\aj+1},R\rangle$ means
$d^\aj\in R^\bot$. Moreover, if $R\notin \Z\cdot R^\ast$, then there is at most
one edge (or a pair of parallel ones) having this property.
\hfill$\Box$\\
\par

%%%%%%%%%%
% (3G.4)
%%%%%%%%

\neu{3G-4}
{\bf Example:}
A typical example of a non-isolated, three-dimensional toric Gorenstein singularity
is the cone over the weighted projective space $\PP(1,2,3)$. We will use it to
demonstrate our calculations of $T^1$ as well as the upcoming construction of
genuine one-parameter families.
$P$ has the vertices $(-1,-1)$, $(2,-1)$, $(-1,1)$, i.e.\ $\sigma$ is generated
from
\[
a^1 =(-1,-1;1)\,,\quad  a^2=(2,-1;1)\,,\quad  a^3=(-1,1;1)\,.
\]
Since our singularity is a cone over a projective variety, $\sigma^{\veee}$ appears
as a cone over some lattice polygon, too. Actually, in our example, $\sigma$ and
$\sigma^{\veee}$ are even isomorphic. We obtain
\[
\sigma^{\veee}=\langle s^1,s^2,s^3\rangle
\quad \mbox{with}\quad
s^1=[0,1;1]\,,\; s^2=[-2,-3;1]\,,\; s^3=[1,0;1]\,.
\]
The Hilbert basis $E\subseteq \sigma^{\veee}\cap\Z^3$ consists of these three
fundamental generators together with
\[
R^\ast=[0,0;1]\,, \quad v^1=[-1,-2;1]\,,\quad v^2=[0,-1;1]\,,\quad w=[-1,-1;1]\,.
\vspace{1ex}
\]
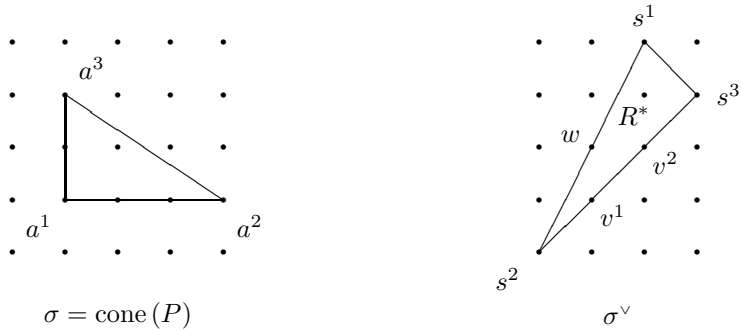
\begin{center}
\unitlength=0.7mm
\linethickness{0.4pt}
\begin{picture}(146.00,65.00)
\put(10.00,20.00){\circle*{1.00}}
\put(10.00,30.00){\circle*{1.00}}
\put(10.00,40.00){\circle*{1.00}}
\put(10.00,50.00){\circle*{1.00}}
\put(10.00,60.00){\circle*{1.00}}
\put(20.00,20.00){\circle*{1.00}}
\put(20.00,30.00){\circle*{1.00}}
\put(20.00,40.00){\circle*{1.00}}
\put(20.00,50.00){\circle*{1.00}}
\put(20.00,60.00){\circle*{1.00}}
\put(30.00,20.00){\circle*{1.00}}
\put(30.00,30.00){\circle*{1.00}}
\put(30.00,40.00){\circle*{1.00}}
\put(30.00,50.00){\circle*{1.00}}
\put(30.00,60.00){\circle*{1.00}}
\put(40.00,20.00){\circle*{1.00}}
\put(40.00,30.00){\circle*{1.00}}
\put(40.00,40.00){\circle*{1.00}}
\put(40.00,50.00){\circle*{1.00}}
\put(40.00,60.00){\circle*{1.00}}
\put(50.00,20.00){\circle*{1.00}}
\put(50.00,30.00){\circle*{1.00}}
\put(50.00,40.00){\circle*{1.00}}
\put(50.00,50.00){\circle*{1.00}}
\put(50.00,60.00){\circle*{1.00}}
\put(110.00,20.00){\circle*{1.00}}
\put(110.00,30.00){\circle*{1.00}}
\put(110.00,40.00){\circle*{1.00}}
\put(110.00,50.00){\circle*{1.00}}
\put(110.00,60.00){\circle*{1.00}}
\put(120.00,20.00){\circle*{1.00}}
\put(120.00,30.00){\circle*{1.00}}
\put(120.00,40.00){\circle*{1.00}}
\put(120.00,50.00){\circle*{1.00}}
\put(120.00,60.00){\circle*{1.00}}
\put(130.00,20.00){\circle*{1.00}}
\put(130.00,30.00){\circle*{1.00}}
\put(130.00,40.00){\circle*{1.00}}
\put(130.00,50.00){\circle*{1.00}}
\put(130.00,60.00){\circle*{1.00}}
\put(140.00,20.00){\circle*{1.00}}
\put(140.00,30.00){\circle*{1.00}}
\put(140.00,40.00){\circle*{1.00}}
\put(140.00,50.00){\circle*{1.00}}
\put(140.00,60.00){\circle*{1.00}}
\put(20.00,50.00){\line(0,-1){20.00}}
\put(20.00,30.00){\line(1,0){30.00}}
\put(50.00,30.00){\line(-3,2){30.00}}
\put(110.00,20.00){\line(1,2){20.00}}
\put(130.00,60.00){\line(1,-1){10.00}}
\put(140.00,50.00){\line(-1,-1){30.00}}
\put(15.00,25.00){\makebox(0,0)[cc]{$a^1$}}
\put(55.00,25.00){\makebox(0,0)[cc]{$a^2$}}
\put(25.00,55.00){\makebox(0,0)[cc]{$a^3$}}
\put(104.00,15.00){\makebox(0,0)[cc]{$s^2$}}
\put(146.00,50.00){\makebox(0,0)[cc]{$s^3$}}
\put(130.00,65.00){\makebox(0,0)[cc]{$s^1$}}
\put(116.00,42.00){\makebox(0,0)[cc]{$w$}}
\put(124.00,27.00){\makebox(0,0)[cc]{$v^1$}}
\put(134.00,37.00){\makebox(0,0)[cc]{$v^2$}}
\put(128.00,46.00){\makebox(0,0)[cc]{$R^\ast$}}
\put(30.00,8.00){\makebox(0,0)[cc]{$\sigma=\mbox{cone}\,(P)$}}
\put(125.00,8.00){\makebox(0,0)[cc]{$\sigma^{\veee}$}}
\end{picture}
\vspace{-2ex}
\end{center}
In particular, $Y_\sigma$ has embedding dimension $7$.
The edges of $P$ have length $\ell(1)=3$, $\ell(2)=1$, and $\ell(3)=2$.
Hence, $Y_\sigma$ contains one-dimensional singularities of transversal type
A$_2$ and A$_1$.\\
According to the previous theorem, $Y_\sigma$ admits only
infinitesimal deformations of the third type. Their degrees come in three series:
\begin{itemize}
\item[($\alpha$)]
% $R^\alpha_p:=
$2R^\ast-p_\alpha\,s^3$ with $p_\alpha\geq 1$. Even the leading
element $R^\alpha=[-1,0,1]$ is not contained in~$\sigma^{\veee}$.
\item[($\beta$)]
% $R^\beta_p:=
$2R^\ast-p_\beta\,s^1$ with $p_\beta\geq 1$. The leading element
equals $R^\beta=v^2=[0,-1,1]$ and sits on the boundary of $\sigma^{\veee}$.
\item[($\gamma$)]
% $R^\gamma_p:=
$3R^\ast-p_\gamma\,s^1$ with $p_\gamma\geq 2$. The leading
element is $R^\gamma=[0,-2,1]\notin\sigma^{\veee}$.
\end{itemize}
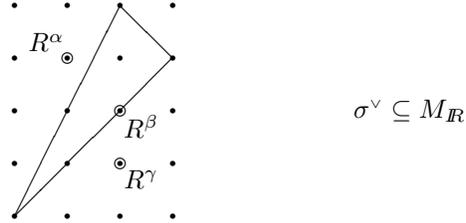
\begin{center}
\unitlength=0.70mm
\linethickness{0.4pt}
\begin{picture}(80.00,50.50)
\put(15.00,10.00){\circle*{1.00}}
\put(15.00,20.00){\circle*{1.00}}
\put(15.00,30.00){\circle*{1.00}}
\put(15.00,40.00){\circle*{1.00}}
\put(15.00,50.00){\circle*{1.00}}
\put(25.00,10.00){\circle*{1.00}}
\put(25.00,20.00){\circle*{1.00}}
\put(25.00,30.00){\circle*{1.00}}
\put(25.00,40.00){\circle*{1.00}}
\put(25.00,50.00){\circle*{1.00}}
\put(35.00,10.00){\circle*{1.00}}
\put(35.00,20.00){\circle*{1.00}}
\put(35.00,30.00){\circle*{1.00}}
\put(35.00,40.00){\circle*{1.00}}
\put(35.00,50.00){\circle*{1.00}}
\put(45.00,10.00){\circle*{1.00}}
\put(45.00,20.00){\circle*{1.00}}
\put(45.00,30.00){\circle*{1.00}}
\put(45.00,40.00){\circle*{1.00}}
\put(45.00,50.00){\circle*{1.00}}
\put(15.00,10.00){\line(1,2){20.00}}
\put(35.00,50.00){\line(1,-1){10.00}}
\put(45.00,40.00){\line(-1,-1){30.00}}
\put(90.00,30.00){\makebox(0,0)[cc]{$\sigma^{\veee}\subseteq M_{\R}$}}
\put(25.00,40.00){\circle{2.00}}
\put(35.00,30.00){\circle{2.00}}
\put(35.00,20.00){\circle{2.00}}
\put(21.00,43.00){\makebox(0,0)[cc]{$R^\alpha$}}
\put(39.00,27.00){\makebox(0,0)[cc]{$R^\beta$}}
\put(39.00,17.00){\makebox(0,0)[cc]{$R^\gamma$}}
\end{picture}
\vspace{-2ex}
\end{center}
\par

%%%%%%%%%%
% (3G.5)
%%%%%%%%

\neu{3G-5}
Each degree belonging to type (3)
(i.e.\ $R=qR^\ast-p\,s^\aj$ with $2\leq q \leq \ell(\aj)$) provides an
infinitesimal deformation. To show that they are unobstructed by describing
how they
lift to genuine one-parameter deformations should be no problem: Just split
the polygon $Q(R)$ into a Minkowski sum meeting conditions
(i) and (ii) of \zitat{Gd}{2}, then construct $\ttau$, $\tsigma$, and
$(f^1-f)$ as in \zitat{Gd}{2} and \zitat{Gd}{5}. However, we prefer to present
the result for our special case all at once by using new coordinates.\\
Let $P\subseteq \A{R^\ast}{}=\R^2\times\{1\}\subseteq \R^3=N_{\R}$ be a lattice
polygon as
in \zitat{3G}{3}, let $R=qR^\ast-p\,s^\aj$ be as just mentioned. Then
$\sigma, \tau\subseteq N_{\R}$ are the cones over $P$ and $P\cap [R\geq 0]$,
respectively, and the one-parameter family in degree $-R$ is obtained as follows:\\
\par

{\bf Proposition:}
{\em
The cone $\ttau\subseteq N_{\R}\oplus\R=\R^4$ is generated by the elements
\begin{itemize}
\item[(i)]
$(a,0)-\langle a,R\rangle\, (\ku{0},1)$, if $a\in P\cap [R\geq 0]$ runs through
the vertices from the $R^\bot$-line until $a^\aj$,
\item[(ii)]
$(a,0)-\langle a,R\rangle \,(d^\aj/\ell(\aj),1)$, if $a\in P\cap [R\geq 0]$ runs 
from $a^{\aj+1}$ until the $R^\bot$-line again, and
\item[(iii)]
$(\ku{0},1)$ and $(d^\aj/\ell(\aj),1)$.
\end{itemize}
The vector space $N_{\R}$ containing $\sigma$  
sits in $N_{\R}\oplus\R$ as $N_{\R}\times\{0\}$. Via this embedding, one obtains
$\tsigma=\ttau+\sigma$ as usual. The monomials $f$ and $f^1$ are given by their
exponents $[R,0], [R,1]\in M\oplus\Z$, respectively.
}
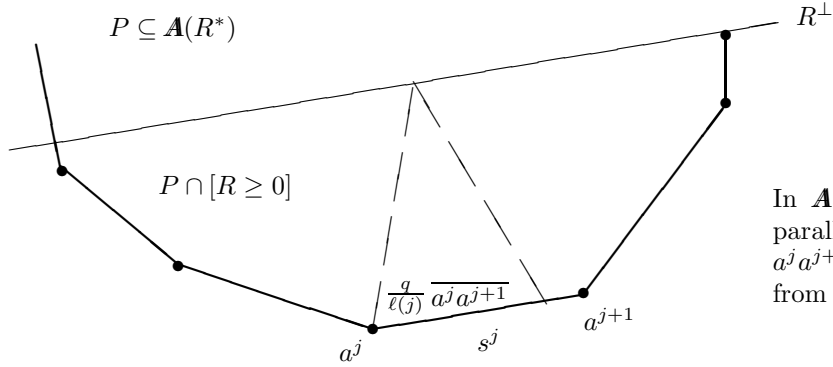
\begin{center}
\unitlength=0.7mm
\linethickness{0.4pt}
\begin{picture}(200.00,66.00)
\thicklines
\put(10.00,60.00){\line(1,-5){4.67}}
\put(14.67,37.00){\line(6,-5){22.33}}
\put(37.00,18.33){\line(3,-1){37.00}}
\put(74.00,6.00){\line(6,1){40.00}}
\put(114.00,12.67){\line(3,4){27.33}}
\put(141.00,49.00){\line(0,1){13.00}}
\put(141.00,62.00){\circle*{2.00}}
\put(141.00,49.00){\circle*{2.00}}
\put(114.00,13.00){\circle*{2.00}}
\put(74.00,6.00){\circle*{2.00}}
\put(37.00,18.00){\circle*{2.00}}
\put(15.00,36.00){\circle*{2.00}}
\put(96.00,4.00){\makebox(0,0)[cc]{$s^\aj$}}
\put(70.00,2.00){\makebox(0,0)[cc]{$a^\aj$}}
\put(119.00,8.00){\makebox(0,0)[cc]{$a^{\aj+1}$}}
\put(158.00,66.00){\makebox(0,0)[cc]{$R^\bot$}}
\put(36.00,63.00){\makebox(0,0)[cc]{$P\subseteq \A{R^\ast}{}$}}
\put(46.00,33.00){\makebox(0,0)[cc]{$P\cap [R\geq 0]$}}
\put(88.00,12.00){\makebox(0,0)[cc]
{\raisebox{0.5ex}{$\frac{q}{\ell(\aj)}$}$\, \ko{a^\aj a^{\aj+1}}$}}
\thinlines
\put(5.00,40.00){\line(6,1){146.00}}
\multiput(74.00,6.00)(2.04,12.24){4}{\line(1,6){1.53}}
\multiput(81.67,53.00)(6.755,-11.258){4}{\line(3,-5){5.066}}
% \put(74.00,6.00){\line(1,6){7.67}}
% \put(81.67,53.00){\line(3,-5){25.33}}
\put(150.00,33.00){\makebox(0,0)[tl]{\parbox{11em}{
   In $\A{R^\ast}{}$, the mutually parallel
   lines $R^\bot$ and $a^\aj a^{\aj+1}$ have distance $q/p$ from each other.}}}
\end{picture}
\end{center}

Geometrically, one can think about $\ttau$ as generated by the interval $I$
with vertices as in (iii) and by the polygon $P^\prime$
obtained as follows: ``Tighten'' $P\cap[R\geq 0]$ along $R^\bot$ by a cone with
base $q/\ell(\aj)\cdot \ko{a^\aj a^{\aj+1}}$ and some top on the $R^\bot$-line;
take $-\langle \bullet, R\rangle$ as an additional, fourth coordinate.
Then, $[R^\ast,0]$ is still $1$ on $P^\prime$ and equals $0$ on $I$.
Moreover, $[R,0]$ vanishes on $I$ and on the $R^\bot$-edge of $P^\prime$;
$[R,1]$ vanishes on the whole $P^\prime$.\\
\par

{\bf Proof:}
We change coordinates. If $g:=\mbox{gcd}(p,q)$ denotes the ``length'' of $R$,
then we can find an $s\in M$ such that $\{s,\,R/g\}$ forms a basis of
$M\cap (d^\aj)^\bot$. Adding some $\kr\in M$ with
$\langle d^\aj/\ell(\aj),\kr\rangle =1$
($\kr$ from Remark \zitat{3G}{3} will do) yields a $\Z$-basis for the whole
lattice $M$. We consider the following commutative diagram:
\[
\dgARROWLENGTH=0.3em
\begin{diagram}
\node{N}
\arrow[4]{e,tb}{(s,\kr,R/g)}{\sim}
\arrow{s,l}{(\mbox{\kf id},\,0)}
\node[4]{\Z^3}
\arrow{s,r}{(\mbox{\kf id},\,g\cdot\mbox{\kf pr}_3)}\\
\node{N\oplus\Z}
\arrow[4]{e,tb}{([s,0],\,[\kr,0],\,[R/g,0],\,[R,1])}{\sim}
\node[4]{\Z^3\oplus\Z}
\end{diagram}
\]
The left hand side contains the data being relevant for our proposition.
Carrying them to the right yields:
\begin{itemize}
\item
$[0,0,g]\in (\Z^3)^\ast$ as the image of $R$;
\item
$[0,0,g,0], [0,0,0,1]\in (\Z^4)^\ast$ as the images of $[R,0]$ and
$[R,1]$, respectively;
\item
$\tau$ becomes a cone with affine cross cut
\vspace{-1ex}
\[
Q([0,0,g])=\mbox{conv}\Big(\big(\langle a,s\rangle/\langle a,R\rangle;\,
\langle a,\kr\rangle/\langle a,R\rangle;\,1/g\big)\,\Big|\;
a\in P\cap [R\geq 0]\Big)\,;
\vspace{-1.2ex}
\]
\item
$I$ changes into the unit interval $(Q_1,1)$ reaching from $(0,0,0,1)$ to
$(0,1,0,1)$;
\item
finally, $\mbox{cone}(P^\prime)$ maps onto the cone spanned by the convex hull
$(Q_0,0)$ of the points
$\big(\langle a,s\rangle/\langle a,R\rangle;\,
\langle a,\kr\rangle/\langle a,R\rangle;\,1/g;\,0\big)$ for $a\in P\cap [R\geq 0]$
on the $a^\aj$-side and\\
$\big(\langle a,s\rangle/\langle a,R\rangle;\,
\langle a,\kr\rangle/\langle a,R\rangle-1;\,1/g;\,0\big)$ for $a$ on the
$a^{\aj+1}$-side, respectively.
\end{itemize}
Since $Q([0,0,g])$ equals the Minkowski sum of the interval
$Q_1\subseteq \A{[0,0,g]}{0}$
and the polygon $Q_0\subseteq\A{[0,0,g]}{}$, we are done by \zitat{Gd}{2}.
\hfill$\Box$\\
\par

%%%%%%%%%%
% (3G.6)
%%%%%%%%

\neu{3G-6}
To see how the original equations of the singularity $Y_\sigma$ will be
perturbed, it is useful to study first the dual cones
$\ttau^{\veee}$ or $\tsigma^{\veee}=\ttau^{\veee}\cap\pi^{-1}(\sigma^{\veee})$:
\vspace{-1ex}\\
\par

{\bf Proposition:}
{\em
If $s\in\sigma^{\veee}\cap M$, then the $(M\oplus\Z)$-element
\[
S:= \left\{ \begin{array}{ll}
[s,\,0] & \mbox{ if } \langle d^\aj,s\rangle\geq 0\\
{}[s,\, -\langle d^\aj/\ell(\aj),\,s\rangle]
& \mbox{ if } \langle d^\aj,s\rangle\leq 0
\end{array} \right.
\]
is a lift of $s$ into $\tsigma^{\veee}\cap (M\otimes\Z)$.
(Notice that it does not depend on $p,q$, but only on $\aj$.)
Moreover, if $s^v$ runs through the edges of $P\cap[R\geq 0]$, the elements
$S^v$ together with $[R,0]$ and $[R,1]$ form the fundamental generators of
$\ttau^{\veee}$.
\vspace{-1ex}
}\\
\par

{\bf Proof:} Since we know $\ttau$ from the previous proposition, the
calculations are straightforward and will be omitted.
\hfill$\Box$\\
\par

%%%%%%%%%%
% (3G.7)
%%%%%%%%

\neu{3G-7}
Recall from \zitat{T1}{1} that $E$ denotes the
minimal set generating the semigroup $\sigma^{\veee}\cap M$.
To any $s\in E$ there is a assigned
variable $z_s$, and $Y_\sigma\subseteq \C^E$ is given by binomial equations
arising from linear relations among elements of $E$.
Everything will be clear by considering an
example: A linear relation such as $s^1+2s^3=s^2+s^4$ transforms into
$z_1\,z_3^2=z_2\,z_4$.\\
The fact that $\sigma$ defines a Gorenstein variety (i.e.\ $\sigma$ is a cone
over a lattice polytope) implies that $E$ consists
only of $R^\ast$ and elements of
$\partial\sigma^{\veee}$ including the fundamental generators $s^v$. If
$E\cap\partial\sigma^{\veee}$ is ordered clockwise, then any two adjacent elements
form together with $R^\ast$ a $\Z$-basis of the three-dimensional lattice $M$.\\
In particular, any three sequenced elements of $E\cap\partial\sigma^{\veee}$
provide a unique linear relation among them and $R^\ast$.
(We met this fact already in Remark \zitat{3G}{3}; there $\kr$, $s^\aj$, and
$\kr^\prime$ were those elements.)
The resulting ``boundary'' equations do not generate the ideal of
$Y_\sigma\subseteq\C^E$. Nevertheless, for describing a deformation
of $Y_\sigma$, it is sufficient to know about perturbations of this subset only.
Moreover, if one has to avoid boundary equations ``overlapping'' a certain spot
on $\partial\sigma^{\veee}$, then it will even be possible to drop up to
two of them from the list.
\vspace{1ex}
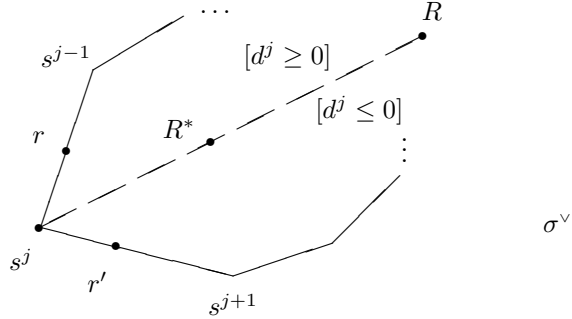
\begin{center}
\unitlength=0.6mm
\linethickness{0.4pt}
\begin{picture}(122.00,68.00)
\put(19.00,55.00){\line(-1,-3){11.67}}
\put(7.33,20.00){\line(4,-1){42.67}}
\put(50.00,9.33){\line(3,1){22.00}}
\put(72.00,16.67){\line(1,1){15.00}}
\put(19.00,55.00){\line(5,3){20.00}}
\put(7.00,20.00){\circle*{2.00}}
\put(13.00,37.00){\circle*{2.00}}
\put(24.00,16.00){\circle*{2.00}}
\put(46.00,68.00){\makebox(0,0)[cc]{$\dots$}}
\put(88.00,39.00){\makebox(0,0)[cc]{$\vdots$}}
\put(38.00,42.00){\makebox(0,0)[cc]{$R^\ast$}}
\put(45.00,39.00){\circle*{2.00}}
\multiput(7.00,20.00)(10.968,5.484){8}{\line(2,1){8.226}}
% \put(7.00,20.00){\line(2,1){85.00}}
\put(92.00,62.50){\circle*{2.00}}
\put(94.00,68.00){\makebox(0,0)[cc]{$R$}}
\put(7.00,40.00){\makebox(0,0)[cc]{$\kr$}}
\put(20.00,8.00){\makebox(0,0)[cc]{$\kr^\prime$}}
\put(3.00,13.00){\makebox(0,0)[cc]{$s^\aj$}}
\put(13.00,58.00){\makebox(0,0)[cc]{$s^{\aj-1}$}}
\put(50.00,3.00){\makebox(0,0)[cc]{$s^{\aj+1}$}}
\put(122.00,21.00){\makebox(0,0)[cc]{$\sigma^{\veee}$}}
\put(62.00,58.00){\makebox(0,0)[cc]{$[d^\aj\geq 0]$}}
\put(78.00,46.00){\makebox(0,0)[cc]{$[d^\aj\leq 0]$}}
\end{picture}
\end{center}

{\bf Theorem:}
{\em
The one-parameter deformation of $Y_\sigma$ in degree $-(q\,R^\ast-p\,s^\aj)$
is completely determined by the following perturbations:
\begin{itemize}
\item[(i)]
(Boundary) equations involving only variables that are induced from
$[d^\aj\geq 0]\subseteq\sigma^{\veee}$ remain unchanged. The same statement holds
for $[d^\aj\leq 0]$.
\item[(ii)]
The boundary equation
$z_\kr\,z_{\kr^\prime}-z_{R^\ast}^{\ell(\aj)}\,z_{s^\aj}^k=0$
assigned to the triple $\{\kr,s^\aj,\kr^\prime\}$
is perturbed
% by the term $t\,z_{R^\ast}^{\ell(\aj)-q}\,z_{s^\aj}^{k+p}$.
into $\big(z_\kr\,z_{\kr^\prime}-z_{R^\ast}^{\ell(\aj)}\,z_{s^\aj}^k\big)
- t\,z_{R^\ast}^{\ell(\aj)-q}\,z_{s^\aj}^{k+p}=0$. Divide everything by
$z^k_{s^\aj}$ if $k<0$.
\vspace{2ex}
\end{itemize}
}
\par

{\bf Proof:}
Restricted to either $[d^\aj\geq 0]$ or $[d^\aj\leq 0]$, the map $s\mapsto S$
lifting $E$-elements into $\tsigma\cap(M\oplus\Z)$ is linear. Hence, any linear
relation remains true, and part (i) is proven.\\
For the second part, we consider the boundary relation
$\kr+\kr^\prime=\ell(\aj)\,R^\ast+k\,s^\aj$ with a suitable $k\in\Z$.
By Lemma \zitat{3G}{6}, the
summands involved lift to the elements $[\kr,0]$, $[\kr^\prime,1]$, $[R^\ast,0]$,
and $[s^\aj,0]$, respectively. In particular, the relation breaks down and has to be
replaced by
\[
\renewcommand{\arraystretch}{1.5}
\begin{array}{rcl}
[\kr,0]+[\kr^\prime,1]&=&
[R,1] + \big(\ell(\aj)-q\big)\, [R^\ast,0] + (k+p)\, [s^\aj,0]
\quad \mbox{ and}\\
\ell(\aj)\,[R^\ast,0]+k\,[s^\aj,0] &=&
[R,0] + \big(\ell(\aj)-q\big)\, [R^\ast,0] + (k+p)\, [s^\aj,0]\,.
\end{array}
\]
The monomials corresponding to $[R,1]$ and $[R,0]$ are $f^1$ and $f$, respectively.
They are {\em not} regular on the total space $X$, but their difference
$t:=f^1-f$ is. Hence, the difference of the monomial versions of both
equations yields the result.\\
Finally, we should remark that (i) and (ii) cover all boundary equations except
those overlapping the intersection of $\partial\sigma^{\veee}$ with
$\ko{R^\ast R}$.
\hfill$\Box$\\
\par

%%%%%%%%%%
% (3G.8)
%%%%%%%%

\neu{3G-8}
We return to Example \zitat{3G}{4} and discuss the one-parameter
deformations occurring in degree $-R^\alpha$, $-R^\beta$, and $-R^\gamma$,
respectively:
\vspace{1ex}\\
{\em Case $\alpha$:}\quad
$R^\alpha=[-1,0,1]=2R^\ast-s^3$ means $\aj=3$, $q=\ell(3)=2$, and $p=1$. Hence,
the line $R^\bot$ has distance $q/p=2$ from its parallel through $a^3$ and
$a^1$. In particular, $\tau=\langle a^1, c^1, c^3, a^3\rangle$ with
$c^1=(1,-1,1)$ and $c^3=(3,-1,3)$.
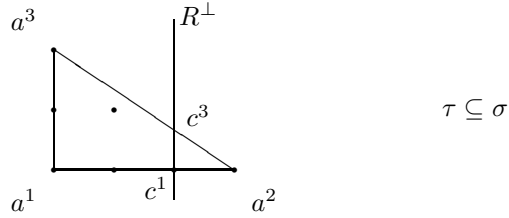
\begin{center}
\unitlength=0.8mm
\linethickness{0.4pt}
\begin{picture}(100.00,47.00)
\put(20.00,15.00){\circle*{1.00}}
\put(20.00,25.00){\circle*{1.00}}
\put(20.00,35.00){\circle*{1.00}}
\put(30.00,15.00){\circle*{1.00}}
\put(30.00,25.00){\circle*{1.00}}
\put(40.00,15.00){\circle*{1.00}}
\put(50.00,15.00){\circle*{1.00}}
\put(20.00,35.00){\line(0,-1){20.00}}
\put(20.00,15.00){\line(1,0){30.00}}
\put(50.00,15.00){\line(-3,2){30.00}}
\put(40.00,10.00){\line(0,1){30.00}}
\put(44.00,41.00){\makebox(0,0)[cc]{$R^\bot$}}
\put(44.00,24.00){\makebox(0,0)[cc]{$c^3$}}
\put(15.00,10.00){\makebox(0,0)[cc]{$a^1$}}
\put(37.00,12.00){\makebox(0,0)[cc]{$c^1$}}
\put(55.00,10.00){\makebox(0,0)[cc]{$a^2$}}
\put(15.00,40.00){\makebox(0,0)[cc]{$a^3$}}
\put(90.00,25.00){\makebox(0,0)[cc]{$\tau\subseteq\sigma$}}
\end{picture}
\vspace{-2ex}
\end{center}
We construct the generators of $\ttau$ by the recipe of Proposition \zitat{3G}{5}:
$a^3$ treated via (i) and $a^1$ treated via (ii) yield the same element
$A:=(-1,1,1,-2)$; from the $R^\bot$-line we obtain $C^1:=(1,-1,1,0)$ and
$C^3:=(3,-1,3,0)$; finally (iii) provides $X:=(0,0,0,1)$ and $Y:=(0,-1,0,1)$.
Hence, $\ttau$ is the cone over the pyramid with plane base $X\,Y\,C^1\,C^3$ 
and $A$ as top. (The relation between the vertices of the quadrangle
equals $3C^1+2X=C^3+2Y$.)
Moreover, $\tsigma$ equals $\tsigma=\ttau+\R_{\geq 0}a^2$ with $a^2:=(a^2,0)$.
Since $A+2X+2a^2=C^3$ and $A+2Y+2a^2=3C^1$, $\tsigma$ is a simplex generated by
$A$, $X$, $Y$, and $a^2$.\\
% \vspace{0.5ex}\\
%
Denoting the variables assigned to $s^1, s^2, s^3, R^\ast, v^1, v^2, w \in E
\subseteq \sigma^{\veee}\cap M$ by $Z_1$, $Z_2$, $Z_3$, $U$, $V_1$, $V_2$, and
$W$, respectively, there are six boundary equations:
\vspace{-1ex}
\[
% fuer \alpha NICHT im Theorem:
Z_3WZ_1-U^3\,=\, Z_1Z_2-W^2\,=\,
% fuer \alpha im Theorem:
WV_1-UZ_2\,=\, Z_2V_2-V_1^2\,=\, V_1Z_3-V_2^2\,=\,
% Teil (ii) vom Theorem:
V_2Z_1-U^2\,=\,0\,.
\vspace{-1ex}
\]
Only the four latter ones are covered by Theorem \zitat{3G}{7}. They will be
perturbed into
\vspace{-1ex}
\[
WV_1-UZ_2 \,=\,Z_2V_2-V_1^2\,=\, V_1Z_3-V_2^2\,=\,
V_2Z_1-U^2-t_\alpha Z_3\,=\,0\,.
\vspace{1ex}
\]
\par

{\em Case $\beta$:}\quad
$R^\beta=[0,-1,1]=2R^\ast-s^1$ means $\aj=1$, $\ell(1)=3$, $q=2$, and $p=1$. 
Hence,
$R^\bot$ still has distance $2$, but now from the line $a^1a^2$.
\vspace{-2ex}
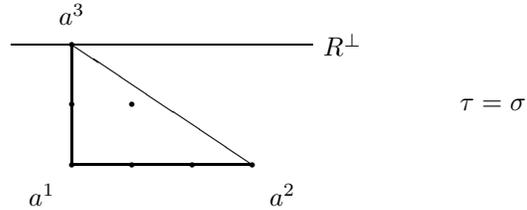
\begin{center}
\unitlength=0.8mm
\linethickness{0.4pt}
\begin{picture}(100.00,47.00)
\put(20.00,15.00){\circle*{1.00}}
\put(20.00,25.00){\circle*{1.00}}
\put(20.00,35.00){\circle*{1.00}}
\put(30.00,15.00){\circle*{1.00}}
\put(30.00,25.00){\circle*{1.00}}
\put(40.00,15.00){\circle*{1.00}}
\put(50.00,15.00){\circle*{1.00}}
\put(20.00,35.00){\line(0,-1){20.00}}
\put(20.00,15.00){\line(1,0){30.00}}
\put(50.00,15.00){\line(-3,2){30.00}}
\put(10.00,35.00){\line(1,0){50.00}}
\put(65.00,35.00){\makebox(0,0)[cc]{$R^\bot$}}
\put(15.00,10.00){\makebox(0,0)[cc]{$a^1$}}
\put(55.00,10.00){\makebox(0,0)[cc]{$a^2$}}
\put(20.00,40.00){\makebox(0,0)[cc]{$a^3$}}
\put(90.00,25.00){\makebox(0,0)[cc]{$\tau=\sigma$}}
\end{picture}
\vspace{-2ex}
\end{center}
We obtain $\ttau=\langle (-1,-1,1,-2); (0,-1,1,-2); (-1,1,1,0);
(0,0,0,1); (1,0,0,1) \rangle$.\\
The boundary equation corresponding to
Theorem \zitat{3G}{7}(ii) is  $Z_3WZ_1-U^3=0$; it perturbs into
$Z_3WZ_1-U^3-t_\beta UZ_1=0$.\\
\par

{\em Case $\gamma$:}\quad
$R^\gamma=[0,-2,1]=3R^\ast-2s^1$ means $\aj=1$, $q=\ell(1)=3$, and $p=2$.
\vspace{-2ex}
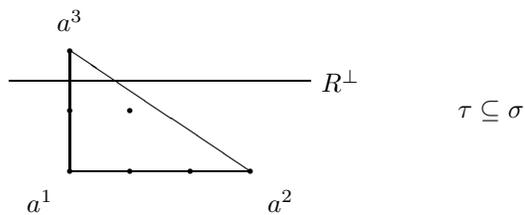
\begin{center}
\unitlength=0.8mm
\linethickness{0.4pt}
\begin{picture}(100.00,47.00)
\put(20.00,15.00){\circle*{1.00}}
\put(20.00,25.00){\circle*{1.00}}
\put(20.00,35.00){\circle*{1.00}}
\put(30.00,15.00){\circle*{1.00}}
\put(30.00,25.00){\circle*{1.00}}
\put(40.00,15.00){\circle*{1.00}}
\put(50.00,15.00){\circle*{1.00}}
\put(20.00,35.00){\line(0,-1){20.00}}
\put(20.00,15.00){\line(1,0){30.00}}
\put(50.00,15.00){\line(-3,2){30.00}}
\put(10.00,30.00){\line(1,0){50.00}}
\put(65.00,30.00){\makebox(0,0)[cc]{$R^\bot$}}
\put(15.00,10.00){\makebox(0,0)[cc]{$a^1$}}
\put(55.00,10.00){\makebox(0,0)[cc]{$a^2$}}
\put(20.00,40.00){\makebox(0,0)[cc]{$a^3$}}
\put(90.00,25.00){\makebox(0,0)[cc]{$\tau\subseteq\sigma$}}
\end{picture}
\vspace{-2ex}
\end{center}
Here, we have $\ttau=\langle (-1,-1,1,-3); (-2,1,2,0); (-1,2,4,0);
(0,0,0,1); (1,0,0,1) \rangle$, and the previous boundary equation provides
$Z_3WZ_1-U^3-t_\gamma Z_1^2=0$.\\
\par

%%%%%%%%%%%%%%%
%
% BIBLIOGRAPHY
%
%%%%%%%%%%%%%%%%%%

%\newpage

\end{document}